\documentclass[11pt]{amsart}

\usepackage{amssymb,latexsym}
\newtheorem{thm}{Theorem}[section]
\newtheorem{lem}[thm]{Lemma}
\newtheorem{defn}[thm]{Definition}

\newtheorem{rem}[thm]{Remark}

\newcommand\R{{\mathbb R}}

\newcommand\m{{\arrowvert}}
\newcommand\n{{\Arrowvert}}
\newcommand\essinf{\mbox{\sl{essinf~}}}
\newcommand\ind{{\mathbf 1}}

\begin{document}

\title{Poisson overlapping microballs: self-similarity and X-ray images}


\author{Hermine BIERM\'E}
\address{MAP5-UMR 8145, Universit\'e Ren\'e Descartes\\
45, rue des Saints-P\`eres  F 75270 PARIS cedex 06  FRANCE \\
hermine.bierme@univ-orleans.fr}

\author{Anne ESTRADE}
\address{MAP5-UMR 8145, Universit\'e Ren\'e Descartes\\
45, rue des Saints-P\`eres  F 75270 PARIS cedex 06  FRANCE\\
anne.estrade@univ-paris5.fr}

\keywords{Random field, random set, overlapping spheres, Poisson point
 process, X-ray transform, asymptotic self-similarity, fractional
 Brownian motion.}
\subjclass{60G60, 60D05, 52A22, 44A12, 60G12, 60G55, 60G57, 60G18, 60H05}.
\maketitle

\begin{abstract}
We study a random field obtained by counting the number of balls
containing each point, when overlapping balls are thrown at random
according to a Poisson random measure. We are particularly
interested in the local asymptotical self-similarity (lass)
properties of the field, as well as the action of X-ray transforms.
We discover two different lass properties when  considering the
asymptotic either {\sl in law} or {\sl on the second order moment}
and prove a relationship between the lass behavior of the field and
the  lass behavior  of its  X-ray transform.  We also  describe  a microscopic
process which leads to a multifractional behavior. These results can be used
to model and analyze porous media, images or connection networks.
\end{abstract}

\section{Introduction}

The purpose of this paper is the  study of a random field obtained by throwing
overlapping balls. Such a field  is particularly well-adapted for modeling
 3D porous or heterogeneous media. In fact we consider
a collection of balls in $\R^3$,  whose centers and radii are chosen  at random according to a Poisson
random measure on $\R^3\times \R^+$. Equivalently, we consider a
germ-grain model where the germs are Poisson distributed and the
grains are balls of random radius.

 The field under study is the mass density defined as the number of balls containing each
point:  the more one point is covered by balls, the higher  is the
mass density at  this point.  From a mathematical point of view, the
dimension three does not yield any specific behavior, so the study
will be carried out in dimension $d$ ($d\ge 1$). Moreover for $d=2$,
the number of balls covering each point defines the discretized gray
level of each pixel in a black and white picture. A one dimensional
($d=1$) germ-grain model is also relevant for modeling communication
networks:  the  germs  stand  for  the  starting  time  of  the
individual ON periods (calls) and the  grains stand  for  the
`half-ball' intervals of duration. The obtained process is a counter
which delivers at each time the number of active connections in the
network.

We have  in mind a microscopic  model which yields  to self-similar
macroscopic properties. In order to get this scaling behavior, we
introduce some power law behavior in the radius distribution and
consider Poisson random measures on $\R^d\times \R^+$ with intensity
of the following type
\begin{equation} \label{intensite}
d\nu(\xi,r)=C r^{-d+\theta}drd\xi~
\end{equation}
for some $\theta$, which may depend on the location $\xi$. The
origin of the `microballs'  described in this paper can be found in
the `micropulses' introduced by Ciosek-Georges and Mandelbrot
\cite{Mand1} with  a fixed power $\theta$ in the intensity measure.
The idea is not new and appeared eighty years ago  when S.D.
Wicksell \cite{Wick} introduced a first model, the famous
`corpuscles', made of random 3D spheres defined as above. The aim of
his study was a stereological question. Since then, this kind  of
model  has been extensively deepened and  extended. We  address to
\cite{SKM} or \cite{Serra} for many examples of random models based
on Poisson point  process and germ-grain models. Let us also mention
two recent papers dealing with similar questions. A one dimensional
germ-grain model  with locations (arrival  times)   uniformly
distributed on the time axis and intervals lengths (call durations) given
by a power law is considered by Cohen and Taqqu in \cite{CT}. A
mixed moving average is performed that sums the height of
connections and the so-called Poissonized Telecom process is
obtained. Also  similar is the model recently studied by Kaj et al.
\cite{KAJ}: the germs are uniformly chosen at random in $\R^d$ and
the grains are obtained by random dilation of a fixed bounded set.
In contrast to the quoted models, let us point out that
inhomogeneity is allowed in our model by choosing a non stationary
intensity measure (\ref{intensite}) with a non constant power
$\theta=\theta(\xi)$.

This paper is  not only concerned with the presentation of  a model
for random media. We also investigate the analysis of the random media mass
intensity following two ways. On one hand, self-similarity
properties are explored. More precisely, we focus on a parameter
that is supposed to contain tangible information on the structure of
the media, the {\it local asymptotical self-similar index}, {\it
lass index} in short. On the other hand, the action of an X-ray
transform on the field is explored. This transform is the
mathematical interpretation for a  radiographic process.  These
techniques  are inspired  from the ones created for Gaussian fields
and are still valuable in the Poisson context. More specifically we
turn to \cite{ABAE} where anisotropic Gaussian fields are analyzed
by performing X-ray transforms and evaluating lass indices.  The
fundamental aim of these methods is to make a 3D  parameter directly
tractable from X-ray images of the
media. \\ \\

The notion of local asymptotic self-similarity has been introduced
in \cite{BJR} in a Gaussian context and extended to the non-Gaussian
realm in  \cite{CL}, or \cite{HB} where a general
presentation is performed for fields with stationary increments. The
lass index can also be related to other parameters of interest as
roughness index \cite{BCI} or Hausdorff dimension \cite{Benassi}. In
the area of network modeling as well, the notion of self-similarity,
at small or large scales, is  fundamental  and  it is highly
connected to long-range dependence. The usual self-similarity
property requires a scale invariance valid for all scales. This is
quite restrictive and we will deal with self-similarity properties
that are fulfilled `at  small  scales' only.  We introduce a light
refinement of the lass property of \cite{BJR}.
\begin{defn}\label{deffddlass}
Let $X=\{X(x); x\in  \R^d\}$ be a random field and $x_0\in \R^d$. We
call the {\it distribution lass index of $X$ at point $x_0$
(fdd-lass index for short)} the following supremum
$$ H_{fdd}(X,x_0)=\sup\left\{\alpha ~;~
\frac{\Delta_{x_0}X(\lambda .)-\mathbf{E}\left(\Delta_{x_0}X(\lambda. )\right)}{\lambda^{\alpha}}\stackrel{fdd}{\rightarrow}
0 ~\mbox{ as }~\lambda\rightarrow 0^+\right\}$$
where   $\Delta_{x_0}X$   denotes  the   field   of   increments  at   $x_0$~:
$$
\Delta_{x_0}X(x)=X(x_0+x)-X(x_0)
$$
 and
$\stackrel{fdd}{\rightarrow}$ means the convergence of the finite dimensional distributions.
\end{defn}
Moreover, when $H=H_{fdd}(X,x_0)$ is finite and when the finite
dimensional distributions of the centered and renormalized
increments $\lambda^{-H}[ \Delta_{x_0}X(\lambda
.)-\mathbf{E}(\Delta_{x_0}X(\lambda. ))]$ converge to the finite
dimensional distributions of a non-vanishing field as
$\lambda\rightarrow 0^+$, the limit field is called {\it tangent
field at point $x_0$} (see \cite{Falc}).

When one deals with  real data, it is almost impossible to  see
whether such a limit exists in distribution. Therefore we introduce
another asymptotic self-similarity property, which only uses the
second  order  moment.
\begin{defn}\label{defcovlass}
Let $X=\{X(x); x\in  \R^d\}$ be a random field and $x_0\in\R^d$. We
call the {\it covariance lass index of $X$ at point $x_0$ (cov-lass
index for short)} the following supremum
$$ H_{cov}(X,x_0)=\sup\left\{\alpha ~;~\mbox{Cov}\left(
\frac{\Delta_{x_0}X(\lambda x)}{\lambda^{\alpha}}, \frac{\Delta_{x_0}X(\lambda x')}{\lambda^{\alpha}}\right)\underset{\lambda\rightarrow 0^+}{\longrightarrow}
0 ~\forall x, x' \in\R^d~\right\}.$$
\end{defn}
By analogy with the fdd situation, when the covariance function of
$\lambda^{-H} \Delta_{x_0}X(\lambda .)$ for $H=H_{cov}$ converges to
a non vanishing covariance function  as $\lambda\rightarrow 0^+$,
the limit covariance will be called the {\it tangent covariance at
point $x_0$}.

Note that the  above self-similarity indices are equal for
 Gaussian fields but not in a general setting. Also note that the existence of a tangent covariance does not imply
existence of the tangent field, and neither  the converse
implication.
 Actually if $H_{cov}$  is the
cov-lass index for $X$ at point $x_0$, then for all $H <H_{cov}$
the covariance function of $\lambda^{-H} \Delta_{x_0}X(\lambda .)$
converges to $0$ as
$\lambda\rightarrow 0^+$. Thus, the finite dimensional distributions of its centered version converges also to $0$ as
$\lambda\rightarrow 0^+$, and
the
fdd-lass index $H_{fdd}$ for $X$ at point $x_0$  - if  exists -  satisfy $H_{fdd}\ge H_{cov}$.

\bigskip

Our main results can be summarized as follows:\\
-  the proposed models provide microscopic descriptions of
macroscopic asymptotical self-similar fields which look like
fractional -or multifractional-
Brownian motions, depending on the involved intensity measure;\\
-  in contrast to the  Gaussian  case, the covariance lass index and
distribution lass index are not equal: the first one can be finite
whereas the second one is infinite, or they can both be finite but
with
different values;\\
-  we present explicit formulas that link the lass indices of a
field and  the  lass  indices of  the  X-ray  transform of  it.  In  particular,  when
  inhomogeneity or anisotropy is introduced in the model, it can be recovered through the
lass indices.

The paper is organized as follows:  the microball model, i.e. the
field that counts the number of balls covering each point, is
introduced in Section 2. The intensities of the Poisson random
measures  we will  use are  specified in  Section \ref{ssmicrob}. A
constant power $\theta$ in (\ref{intensite}) will yield to a
fractional microball model which is stationary and isotropic. A
non-constant power $\theta$ will yield the multifractional model.
 We also introduce  in Section \ref{ssXray} the X-ray transform. Section 3 is devoted to the self-similarity
properties  of  the  microball  model  and its  X-ray  transform.
Theorem \ref{thlassm} deals  with the fractional model and  Theorems
\ref{thlassh} and \ref{thlasshstar}    deal    with     the
multifractional    model,    where $\theta=\theta(\xi)$  is a
smooth, or a singular  function respectively. We  also compare our
results to homogenization results in Section \ref{sshom}. The
proofs of Theorems \ref{thlassh} and \ref{thlasshstar} are detailed
at the end of the paper.  The covariance  lass properties  are
proved  in Section \ref{covproof},  and the distribution lass properties in
Section \ref{fddproof}.

\section{ The microball model and its X-ray transforms}
\subsection{The random grain model}

The model is built by considering the superposition of balls $B(\xi,r)$, where
$\xi$ is a point in $\R^d$ and $r>0$ is the radius of the ball.
As in \cite{KAJ}, we want to study the mass distribution generated by a family
of balls $B(\xi_j,r_j)$, with random location $\xi_j$ and random radius $r_j$. We assume that $(\xi_j,r_j)$ are given by a Poisson point process with intensity
$\nu(d\xi,dr)$, where $\nu$ is a non-negative $\sigma$-finite measure on $\R^d\times\R^+$.

For each $x\in\R^d$  we are naturally interested in the number of
balls $B(\xi,r)$ that contain the point $x$, given by
$$
\#\left\{j; x\in B(\xi_j,r_j)\right\}
=\sum_j\mathbf{1}_{B(\xi_j,r_j)}(x).
$$
Such a field is well defined as soon as
\begin{equation} \label{int}
\forall x \in \R^d~,~ \int_{\R^d\times
  \R^+}\mathbf{1}_{B(\xi,r)}(x)\nu(d\xi,dr)~ <~ +\infty~.
\end{equation}
Moreover, in that case, we can represent the field through a
stochastic integral with respect to a Poisson measure $N$ with
intensity $\nu$, as
$$
 \int_{\R^d\times
  \R^+}\mathbf{1}_{B(\xi,r)}(x)N(d\xi,dr).
$$
Let  us also  describe  the intuitive  scenario we  have  in mind
for the one dimensional case ($d=1$). We look at the
`half-ball'-interval $[\xi,\xi+r)$ as the ON period of a single
call and the number of connected users at time $x$ is equal to the
one dimensional integral -if exists- $\int_{\R\times
\R^+}\mathbf{1}_{[\xi,\xi+r)}(x)N(d\xi,dr)$.  Since the integrand
$\mathbf{1}_{[\xi,\xi+r)}$   behaves  like $\mathbf{1}_{B(\xi,r)}$,
we will not distinguish anymore the special case $d=1$.

\subsection{ The microball model} \label{ssmicrob}
 As  in \cite{Mand1},  we
assume that the radii of such a random grain model  follow a power
law in $(0,1)$.
 The exponent of the power law can either be
constant (fractional case) or can depend on the center of the ball
$\xi$ (multifractional case). More precisely we consider intensity
measures on $\R^d\times \R^+$ with special shapes $\nu_m$ and
$\nu_h$ described below.

For  $m>0$  we define the  fractional  intensity  measure $\nu_m$ on
$\R^d\times \R^+$ as
\begin{equation} \label{nu_m}
\nu_{m}(d\xi,dr)=r^{-d-1+2m}\mathbf{1}_{(0,1)}(r)d\xi dr~.
\end{equation}
For  a  function $h$ on $\R^d$ such that ${\essinf}h~ >0,$ we define
the multifractional intensity measure $\nu_h$ on $\R^d\times \R^+$
as
\begin{equation} \label{nu_h}
\nu_{h}(d\xi,dr)=r^{-d-1+2h(\xi)}\mathbf{1}_{(0,1)}(r)d\xi dr~.
\end{equation}
It is straightforward to see that these kind of measures satisfy
(\ref{int}).
 This allows the following definition.
\begin{defn} \label{def}
Let $h$ be a function  on $\R^d$ such that ${\essinf}h~ >0,$  and
let ${N_h}$ be a Poisson random measure with intensity $\nu_h$. The
field defined on $\R^d$ as
$$
X(x)=\int_{\R^d\times
  \R^+}\mathbf{1}_{B(\xi,r)}(x){N}_h(d\xi,dr)$$
is called a  {\it microball model with index $h$. }
\end{defn}

Note that the microball model  has moments of all order. In particular, its mean value is given, for each $x\in\R^d$, by
$$\mathbf{E}(X(x))=\int_{\R^d\times \R^+}\mathbf{1}_{B(\xi,r)}(x)\nu_h(d\xi,dr).$$
Moreover, by the isometry of the Poisson measure, its covariance
function is equal to
$$\mbox{Cov}(X(x),X(x'))=\int_{\R^d\times \R^+}\mathbf{1}_{B(\xi,r)}(x)\mathbf{1}_{B(\xi,r)}(x')\nu_h(d\xi,dr),$$
for all $x,~x'\in \R^d$.

\bigskip

The self-similarity properties of the microball model that we will
study in the sequel deal with the local behavior of $X$. We now
compute the increments of $X$ and analyse their moments in the
following lemma. In what follows we write for $x_0, x\in\R^d$
\begin{equation}\label{deltaX}
\Delta_{x_0}X(x)= X(x_0+x)- X(x_0)=
\int_{\R^d\times\R^+}\psi(x,\xi-x_0,r){N}_h(d\xi,dr),
\end{equation}
where
\begin{equation}\label{psi}
\psi(x,\xi,r)=\mathbf{1}_{B(\xi,r)}(x)-\mathbf{1}_{B(\xi,r)}(0)
=\mathbf{1}_{\m x-\xi\m<r\le \m \xi\m}-\mathbf{1}_{ \m \xi\m <r \le \m x-\xi\m},
\end{equation}
and $\m.\m$ denotes the usual Euclidean norm.

\begin{lem}
\label{moment de psi} Let $m\in (0,1/2)$. There exists a constant
$C(m)\in (0,+\infty)$ such that, for all $p\in (0,+\infty)$  and
$x\in \R^d$,
$$\int_{\R^d \times \R^+}\m \psi(x,\xi,r)\m^p ~r^{-d-1+2m}dr d\xi=C(m)\m x\m^{2m}.$$
\end{lem}
\begin{proof}[Proof] Note that $\m \psi(x,\xi,r)\m^p=\m \psi(x,\xi,r)\m=\mathbf{1}_{\m x-\xi\m<
r\le\m\xi\m}+\mathbf{1}_{\m \xi\m< r\le \m x-\xi\m}.$\\
Hence, for all $x, \xi\in \R^d$,
$$\int_{\R^+} \m \psi(x,\xi,r)\m^p ~ r^{-d-1+2m} dr=\frac{2}{d-2m}\mathbf{1}_{\m \xi\m<  \m x-\xi\m}
\left(\m \xi\m^{-d+2m}-\m x-\xi\m^{-d+2m}\right)$$
and for  $0<m<1/2$, the function $\xi \mapsto \m \xi\m^{-d+2m}-\m x-\xi\m^{-d+2m}$ is integrable on $\R^d$.
The conclusion is obtained by rotation invariance and homogeneity.
\end{proof}

\subsection{X-Ray transform} \label{ssXray}

One motivation for this paper is to describe, model and analyze 
 heterogeneous  media.  We  have  in   mind  the  possibility  to  estimate  a
 macroscopic 3D parameter through X-ray images. By this method, it will be possible
to get an analysis of the media without entering  the media
(non-invasive method). In this section the mathematical tool
associated with X-ray images is presented and tested on the
microball model. We assume that $d\ge 2$.

Following  the usual  notation  (see \cite{ramm}  for  instance),
for a direction $\alpha\in S^{d-1}=\left\{x\in\R^d; \m
x\m=1\right\}$, the X-ray transform in the direction $\alpha$ of any
function $f\in L^1(\R^d)$ is given by
$$y\in <\alpha>^{\perp}~ \mapsto ~ \int_{\R}f(y+p\alpha)dp$$
where $ <\alpha>^{\perp}:=\left\{x\in\R^d; x\cdot\alpha=0\right\},$
and $\cdot$ denotes the usual scalar product on $\R^d$. We want to
define, in the same way, the X-ray transform of a microball model
$X$. Unfortunately, the realizations $x \in \R^d \mapsto
X(x,\omega)$ do not belong to $L^1(\R^d)$. We will therefore work
with the {\it windowed } X-ray transform defined through a fixed
window $\rho$. We assume that $\rho$ is a continuous function on
$\R$ with fast decay, namely
\begin{equation}\label{decroissance rapide}
\forall N\in{\mathbb N},\,\exists C_N,\, \forall p\in \R,\,\,\, \left|\rho(p)\right|\le C_N\left(1+\m p\m\right)^{-N}.
\end{equation}
 For any function $f\in L^1_{loc}(\R^d)$ with slow growth, we define the {\it windowed X-ray transform of $f$ in the direction
$\alpha$} to be the map
\begin{equation} \label{eqPf}
y\in <\alpha>^{\perp}~\mapsto ~{\mathcal
P}_{\alpha}f(y):=\int_{\R}f(y+p\alpha)\rho(p)dp~.
\end{equation}

Our aim is to apply such a transformation to a microball model.
We use the properties of the Poisson random measure and  the fact
that one can define the field
$$\left\{\int_{\R^d\times\R^+}\left(\int_{\R}\mathbf{1}_{B(\xi,r)}(y+p\alpha)\rho(p)dp\right)N_h(d\xi,dr); ~y~\in <\alpha>^{\perp}~\right\},$$
as soon as the integrand
$${\mathcal P}_{\alpha}\left(\mathbf{1}_{B(\xi,r)}(.)\right)(y)=\int_{\R}\mathbf{1}_{B(\xi,r)}(y+p\alpha)\rho(p)dp$$
is integrable with respect to $\nu_h(d\xi,dr)$ for each $y\in <\alpha>^{\perp}$. But this is straightforward using (\ref{int}) and (\ref{decroissance rapide}). This allows us to state the following definition.
\begin{defn}\label{defXray}
Let $h$ be a function  on $\R^d$ such that ${\essinf} ~h~ >0$ and
let ${N_h}$ be a Poisson random measure with intensity $\nu_h$. Let
$X$ be the microball model with index $h$,
$X(\cdot)=\int_{\R^d\times
  \R^+}\mathbf{1}_{B(\xi,r)}(\cdot){N}_h(d\xi,dr).$\\
Let $\rho$ be a continuous window that satisfies (\ref{decroissance
rapide}). For $\alpha\in S^{d-1}$, we call the windowed X-ray
transform of $X$ in the direction $\alpha$, the field defined on
$<\alpha>^{\perp}$ by
$${\mathcal P}_{\alpha}X(y)=\int_{\R^d\times \R^+}{\mathcal P}_{\alpha}\left(\mathbf{1}_{B(\xi,r)}(\cdot)\right)(y){N}_h(d\xi,dr).$$
\end{defn}

Note that for $(\xi,r)$ in $\R^d\times \R^+$, Cauchy-Schwarz inequality leads to
\begin{eqnarray*}
\left|{\mathcal
P}_{\alpha}\left(\mathbf{1}_{B(\xi,r)}(\cdot)\right)(y) \right|^2
&\le & 2 (r^2-\m
y-\pi_{\alpha^{\perp}}(\xi)\m^2)^{1/2}_+~\int_{\R}\mathbf{1}_{B(\xi,r)}\left(y+p\alpha\right)\rho^2(p)dp \\
&\le& 2 r\int_{\R}\mathbf{1}_{B(\xi,r)}\left(y+p\alpha\right)\rho^2(p)dp~,
\end{eqnarray*}
where
$\pi_{\alpha^{\perp}}$ denotes the orthogonal projection on
$<\alpha>^{\perp}$. Thus, for all $y\in <\alpha>^{\perp}$, the function
${\mathcal   P}_{\alpha}\left(\mathbf{1}_{B(\xi,r)}(.)\right)(y)$  belongs  to
$L^2\left(\R^d\times\R^+,\nu_h(d\xi,dr\right)$
and ${\mathcal P}_{\alpha}X$ admits a second order moment. \\

\bigskip
As for the microball model itself, we are particularly interested in
the local behavior of the X-ray transform and we need to estimate
its increments.

For $y_0, y\in <\alpha>^{\perp}$, using (\ref{deltaX}) we get
\begin{equation} \label{deltaP}
\Delta_{y_0}{\mathcal P}_{\alpha}X(y)={\mathcal P}_{\alpha}X(y_0+y)-{\mathcal P}_{\alpha}X(y_0)
=\int_{\R^d\times \R^+}
G_{\rho}(y,\xi-y_0,r){N}_h(d\xi,dr)
\end{equation}
where
\begin{equation} \label{grho}
G_{\rho}(y,\xi,r)= \int_{\R}\psi(y,\xi-p\alpha,r)\rho(p)dp~.
\end{equation}
Note that the above integral is well defined for any bounded
function $\rho$. In the special case where $\rho \equiv 1$, we write
$G$ instead of $G_1$ and for $y, \gamma$ in $<\alpha>^{\perp}$ and
$r$ in $\R^+$, a simple computation gives
\begin{equation} \label{eqG}
G(y,\gamma,r)=G_1(y,\gamma,r)=(r^2-\m     y-\gamma\m^2)^{1/2}_{+}     -(r^2-\m
\gamma\m^2)^{1/2}_{+},
\end{equation}
where, as usual, $t_+:=\max(0,t)$ for all $t\in \R$.\\
The next lemma provides upper-bounds for the integral of $G(y,.)$.

\begin{lem}\label{GdansL2}
For  $m\in (0,1/2)$ and $y\in <\alpha>^{\perp}$,
$$G(y,.)\in
L^2(<\alpha>^{\perp}\times \R^+,r^{-d-1+2m}d\gamma dr)~.$$
\end{lem}

\begin{proof}[Proof]
For $m\in (0,1/2)$, $y\in <\alpha>^{\perp}$ and $\lambda>0$, on one
hand
\begin{eqnarray}
\int_{<\alpha>^{\perp}}
G(y,\gamma,r)^2   d\gamma
&\le & r^2\int_{<\alpha>^{\perp}} \left(\mathbf{1}_{\m
y-\gamma\m< r}+\mathbf{1}_{\m \gamma\m< r}\right)~
d\gamma \nonumber \\
&\le &C~ r^{d+1}. \label{integrale en 0}
\end{eqnarray}
On the other hand,  a change of variable gives, for  $y\ne 0$ and
$r>0$,
$$
\left(\frac{\m y\m }{2}\right)^{-(d+1)}\int_{<\alpha>^{\perp}}
G(y,\gamma,r)^2 ~d\gamma ~$$
$$=\int_{<\alpha>^{\perp}}    \left|     \left( \left(\frac{2r}{\m y\m}\right)^2-\left|
    \gamma-\frac{y}{\m y\m}\right|^2 \right)_+^{1/2}-
 \left( \left(\frac{2r}{\m y\m}\right)^2-\left| \gamma +\frac{y}{\m y\m}\right|^2 \right)
  _+^{1/2} \right|^2~ d\gamma ~.$$
The next lemma provides an upper bound for the last quantity, which
leads to
\begin{equation} \label{intG2}
\int_{<\alpha>^{\perp}}
G(y,\gamma,r)^2 ~d\gamma \leq C |y|^2 r^{d-1} \ln (2+\frac{2r}{|y|})~.
\end{equation}
Since $m\in (0,1/2)$,  inequalities (\ref{integrale en 0}) and
(\ref{intG2}) conclude for the proof.
\end{proof}
\begin{lem} \label{lemaline}
Let $n\in {\mathbb N}^*$. There exists a constant $C>0$ such that
for all direction $e\in S^{n-1}$ and all $r>0$,
$$\int_{\R^{n}}\left| (r^2-\m x-e\m^2)_+^{1/2}-(r^2-\m x+e\m^2)_+^{1/2}\right|^2dx \le Cr^n\ln(2+r).$$
\end{lem}
\begin{proof}
For $n=1$, we have to prove that there exists a constant $C$ such
that, for $r>0$,
$$\int_0^{r+1}\left\vert (r^2-(x-1)^2)_+^{1/2}-(r^2-(x+1)^2)_+^{1/2}
\right\vert^2dx
\leq Cr \ln (r+2).$$ This is an easy consequence of the fact that the
function  that we integrate is bounded
 by $4r$ for $x\in [r-1,r+1]$, and
by $16r^2((r-1)(r-x+1))^{-1}$ for  $x\in [0,r-1]$ when $r>1$. \\
In   the    general   case   ($n>   1$)   we   write    $x=x'+x''e$
with $x'=\pi_{<e>^{\perp}}(x)$ and $x''\in \R$. From the
one-dimensional case, for $x'\in <e>^{\perp}$~,
$$ \int_{\R}\left| (r^2-\m x'\m^2- \m x''-1\m^2)_+^{1/2}-(r^2
-\m x'\m^2-\m x''+1\m^2)_+^{1/2}\right|^2~dx"$$
$$ \le C(r^2-\m x'\m^2)_+^{1/2}\ln(r+2).$$
But
$$\int_{<e>^{\perp}}(r^2-\m x'\m^2)_+^{1/2}~dx'=r^n\m S^{n-2}\m \int_0^1(1-t^2)^{1/2}t^{n-2}dt.$$
Finally, we can change the constant $C$ such that
$$\int_{\R^{n}}\left| (r^2-\m x-e\m^2)_+^{1/2}-(r^2-\m x+e\m^2)_+^{1/2}\right|^2~dx \le Cr^n\ln(2+r).$$
\end{proof}

\section{Lass properties}

The section is  devoted to the study of  self-similarity properties
of the microball models and their X-ray transforms.

\subsection{The fractional microball model} \label{ssfmm}

Let $d\ge 1$. Let $m>0$ and recall that the   fractional  intensity  measure  $\nu_m$  on
$\R^d\times \R^+$ is given by
$$
\nu_{m}(d\xi,dr)=r^{-d-1+2m}\mathbf{1}_{(0,1)}(r)d\xi dr
$$
and the fractional microball model by
$$X=\left\{\int_{\R^d\times\R^+}\mathbf{1}_{B(\xi,r)}(x)N_m(d\xi,dr);~x~\in\R^d\right\},$$
where $N_m$ is a random Poisson measure of intensity $\nu_m$.

\medskip

First, let  us remark that the choice  of $\nu_m$ as intensity  means that the
 center and the radius are thrown independently.
Hence, since the centers are  uniformly distributed on the state space $\R^d$,
 the fractional microball model $X$ is isotropic and stationary, i.e. for every rotation $R$
centered   at  $0$  in   $\R^d$  and   for  all   $x_0\in \R^d$,
$$X\circ R\stackrel{fdd}{=}X ~\mbox{ and }~ X(x_0+.)\stackrel{fdd}{=}X~.$$

\medskip

In the following, we will assume that $m<1/2$ so that Lemma
\ref{moment de psi} applies. Then, using (\ref{deltaX}) and the
stationarity of $X$ we get for $x, x'\in \R^d$,
$$\mathbf{E}\left(X(x)-X(x')\right)^2\le C(m)\m x-x'\m^{2m}.$$
 Thus, the field $X$ is mean square continuous.

Furthermore, by the correlation theory of stationary random fields (see \cite{Yaglom} for example), there exists
a finite positive Radon measure $\sigma$ such that
$$\mbox{Cov}\left(X(x),X(0)\right)=\int_{\R^d}e^{-ix.\xi}d\sigma(\xi).$$
By computing the inverse Fourier transform we obtain that $\sigma$
is absolutely continuous with respect to the Lebesgue's measure.
Moreover, the spectral density of the fractional microball  model
$X$ is given by
$$\left(2\pi\right)^{d/2}\m \xi\m^{-2m-d}\int_0^{\m\xi\m}J_{d/2}^2(s)s^{-d-1+2m}ds,$$
where $J_{d/2}$ is the Bessel function (see \cite{ramm} p. 406 for example). Let us recall that
$$J_{d/2}(s)=\frac{s^{d/2}}{2^{d/2}\Gamma(1+d/2)}\left(1+{\mathcal O}\left(\m s\m^2\right)\right),\,\, s\rightarrow 0,$$
and
$$J_{d/2}(s)=\sqrt{\frac{2}{\pi s}}\left(\cos\left(s-\frac{(d+1)\pi}{4}\right)+{\mathcal O}\left(s^{-1}\right)\right),\,\, s\rightarrow +\infty.$$
Hence, the previous  spectral density has a power law behavior and lass properties for the associated field are expected.  This point of view is detailed in \cite{HB}
section 2.2.2. \\

We first  investigate the  covariance lass  property of  $X$. Recall
that $X$ is stationary and therefore the behavior of the increments
of $X$ around any $x_0 \in \R^d$ does not depend on $x_0$. We have
to deal, for $x, x' \in\R^d$ and $\lambda\rightarrow 0^+$, with
$$\mbox{Cov}(\Delta_{x_0}X(\lambda x),\Delta_{x_0}X(\lambda x'))=\int_{\R^d\times (0,1)}\psi(\lambda x,\xi,r)\psi(\lambda x',\xi,r)  ~r^{-d-1+2m}dr
d\xi~,$$
where $\psi$ is given by (\ref{psi}).
By homogeneity and integrability of $\psi$ (see Lemma \ref{moment de psi}) it is easy
to establish that, for all $x,x'\in \R^d$,
\begin{eqnarray}
&\lim_{\lambda \rightarrow 0^+} & \lambda^{-2m}\int_{\R^d \times  (0,1)} \psi(\lambda  x,\xi,r)\psi(\lambda x',\xi,r)
~r^{-d-1+2m}dr d\xi  \nonumber\\
&=& \int_{\R^d   \times  \R^+}   \psi(x,\xi,r)\psi(x',\xi,r)  ~r^{-d-1+2m}dr
d\xi~. \label{limpsi}
\end{eqnarray}
The cov-lass  property of  $X$ follows from (\ref{limpsi}).

Now, let us  assume that $d\ge 2$ and choose $\alpha\in  S^{d-1}$. We look for
the cov-lass property of ${\mathcal P}_{\alpha}X$,
the windowed X-ray transform of $X$. Let us recall (\ref{deltaP}) and (\ref{grho}): for any $y_0$ and $y$ in $<\alpha>^{\perp}$,
$$\Delta_{y_0}{\mathcal P}_{\alpha}X(y)
=\int_{\R^d\times \R^+}
G_{\rho}(y,\xi-y_0,r){N}_m(d\xi,dr)$$
where
$$ G_{\rho}(y,\xi,r)= \int_{\R}\psi(y,\xi-p\alpha,r)\rho(p)dp~. $$
It is straightforward to see that ${\mathcal P}_{\alpha}X$ is also stationary so that the cov-lass properties will not depend on the point $y_0$.

Writing any $\xi \in \R^d$ as $\xi=\gamma + t\alpha$ with
$\gamma \in <\alpha>^{\perp}$ and $t\in \R$ and performing a  change
of variable (translation-dilation) in the integral that defines
$G_{\rho}$, we obtain
\begin{equation}\label{covxray}
\mbox{Cov}(\Delta_{y_0}{\mathcal P}_{\alpha}X(\lambda y),\Delta_{y_0}{\mathcal P}_{\alpha}X(\lambda y'))
\end{equation}
$$=\lambda^{1+2m}
\int_{<\alpha>^{\perp}\times \R\times (0,\lambda^{-1})}
G_{\rho(t+\lambda .)}(y,\gamma,r)G_{\rho(t+\lambda .)}(y',\gamma,r)
~r^{-d-1+2m}~d\gamma dt dr, $$
with $\rho(t+\lambda .)$ denoting the window $p\mapsto \rho(t+\lambda p)$.
Note that for  $r\in\R^+$, $t\in \R$ and $y,\gamma \in <\alpha>^{\perp}$,
$$G_{\rho(t+\lambda .)}(y,\gamma,r) ~\stackrel{\lambda\rightarrow 0^+}{\longrightarrow} ~G(y,\gamma,r) \rho(t)~.$$
Assumption (\ref{decroissance rapide}) on $\rho$ and Lemma \ref{GdansL2} allow to conclude that
$$\lim_{\lambda\rightarrow 0^+}\lambda^{-1-2m}\mbox{Cov}(\Delta_{y_0}{\mathcal P}_{\alpha}X(\lambda y),\Delta_{y_0}{\mathcal P}_{\alpha}X(\lambda y'))$$
$$=\left(\int_{\R}\rho(t)^2dt\right)
\int_{<\alpha>^{\perp}\times \R^+}G(y,\gamma,r)G(y',\gamma,r)r^{-d-1+2m}~d\gamma  dr.$$

We now state a theorem including the  notions of local asymptotical
self-similarity    index   given    in    Definitions   \ref{deffddlass}    and
\ref{defcovlass}. The fdd-lass properties are proved in Section \ref{fddproof}, as a consequence of Theorem \ref{thlassh} in the special
case where the multifractional index $h$ is constant equal to $m$.
\begin{thm} \label{thlassm}
 Let $d\ge 2$. Let $X$ be the fractional microball model with index $m \in  \left(0,1/2\right)$
and let  ${\mathcal P}_{\alpha}X$  be its $X$-ray  transform in  the direction
$\alpha \in S^{d-1}$.
\begin{itemize}
\item At any point
$x_0\in \R^d$, the cov-lass index of $X$ is equal to $m$ and the fdd-lass index of $X$ is equal to $+\infty$. Moreover the
covariance  of   $\lambda^{-m}\Delta_{x_0}X(\lambda.)$  converges,  up   to  a
multiplicative constant, to the
covariance of a fractional Brownian motion of index $m$.
\item  At any  point $y_0\in  <\alpha>^{\perp}$,  the cov-lass  index and  the
  fdd-lass index of ${\mathcal P}_{\alpha}X$  are equal to $m+1/2$.  Moreover the
covariance and the finite dimensional distributions of $\lambda^{-m-1/2}\left(\Delta_{y_0}({\mathcal P}_{\alpha}X) (\lambda.)-\mathbf{E}\left(\Delta_{y_0}({\mathcal P}_{\alpha}X) (\lambda.)\right)\right)$ converge, up to a
multiplicative constant, to the corresponding ones  of a fractional Brownian motion of index $m+1/2$.
\end{itemize}

\end{thm}

\begin{rem} The first point of Theorem \ref{thlassm} is still true in the one-dimensional case ($d=1$).
\end{rem}

Let us comment this result.

In  one dimension, this result describes the  small scale behavior
of  the number of active  connections in a communication network:
the covariance is locally asymptotically self similar and behaves
like a fractional Brownian motion covariance. More generally, the
same is observed in the multi-dimensional  case. Hence the
fractional microball model provides a microscopic description of a
random media which behaves, up to the second-order moment, like a
fractional Brownian motion.

The second point of Theorem \ref{thlassm} is very interesting from a
practical point of view. The one-to-one correspondence between the
lass index of $X$ and the lass index of ${\mathcal P}_{\alpha}X$
allows the estimation of the 3D lass index through the analysis of
the media radiographic images. The  fractional microball  model is
thus relevant  to model isotropic and stationary media.

 Following a widespread idea (\cite{Peltier},
\cite{BJR},\cite{ABAE}), we now consider the multifractional
microball model, where  $m$ is replaced by a function that depends
on the ball location.

\subsection{The multifractional microball model} \label{ssmfmm}

  In this section we study the lass properties of the multifractional microball model, associated with  the   multifractional  intensity  measure  $\nu_h$,
given on $\R^d\times \R^+$   by
$$
\nu_{h}(d\xi,dr)=r^{-d-1+2h(\xi)}\mathbf{1}_{(0,1)}(r)d\xi dr~,
$$
where $h$ is a  function  on $\R^d$ such that $0< h(\xi) <1/2~$.
\medskip

First, let  us remark that the multifractional microball model is
not stationary nor isotropic when $h$ is not constant.

Let us  deal with the cov-lass properties of the multifractional microball. One has to study,
for all  $x_0, x,  x' \in\R^d$, and  $\lambda\rightarrow 0^+$,  the asymptotic
behavior of $\mbox{Cov}\left(
{\Delta_{x_0}X(\lambda x)}, {\Delta_{x_0}X(\lambda x')}\right)$. By a change of variables, with
$\psi$ given by (\ref{psi}),
\begin{equation}\label{covmmm}
\mbox{Cov}\left(
{\Delta_{x_0}X(\lambda x)}, {\Delta_{x_0}X(\lambda x')}\right)
\end{equation}
$$=\int_{\R^d\times
(0,\lambda^{-1})}\lambda^{2h(x_0+\lambda\xi)}\psi( x,\xi,r)\psi(
x',\xi,r)r^{-d-1+2h(x_0+\lambda\xi)}d\xi dr.$$
Similarly, when $d\ge 2$ and  $\alpha\in S^{d-1}$, for all $y_0, y, y' \in<\alpha>^{\perp}$,  (\ref{covxray}) becomes
\begin{equation}\label{covxraymmm}
\mbox{Cov}(\Delta_{y_0}{\mathcal P}_{\alpha}X(\lambda y),\Delta_{y_0}{\mathcal
  P}_{\alpha}X(\lambda y'))
\end{equation}
$$=\lambda~
\int_{<\alpha>^{\perp}\times \R\times (0,\lambda^{-1})}
G_{\rho(t+\lambda .)}(y,\gamma,r)G_{\rho(t+\lambda .)}(y',\gamma,r)$$
$$ \times~ (\lambda r)^{2h(y_0+t\alpha+\lambda\gamma)}r^{-d-1}~d\gamma dt dr.$$

In  order to  get cov-lass  properties  for both  $X$ and  its
windowed  X-ray transform, further assumptions on $h$ have to be
made. We are mainly  interested in two kinds of functions. The first
kind deals with the case where $h$ is smooth on $\R^{d}$. It is
linked with the  multifractional Brownian  motion \cite{Peltier},
\cite{BJR}, obtained by substituting the Hurst parameter $H$ by a
Lipschitz function  on the state space. The second class of function
is when $h$ is singular at the point 0,  defined on
$\R^d\smallsetminus\{0\}$ by homogeneity, ie $h(\lambda\xi)=h(\xi)$
for all $\lambda\in\R^*$. This follows the point of
 view of \cite{ABAE} to get anisotropic
generalizations of the fractional Brownian motion.\\

Let us recall the definition of a $\beta$-Lipschitz function.
\begin{defn}
Let $\left(M,d_M\right)$ be a metric space and $\beta \in (0,1]$. A
function $f:M\rightarrow \R$ is called $\beta$-Lipschitz on $M$ if
there exists $C>0$ such that
$$\forall x, y\in M,\,\,\,\,d_M(x,y)\le 1\Rightarrow \left|f(x)-f(y)\right|\le Cd_M(x,y)^{\beta}.$$
\end{defn}

\subsubsection{The smooth case:}

let us first study the case of a $\beta$-Lipschitz function on
$\R^d$. By continuity of $h$ around $x_0\in\R^d$,
we get from (\ref{covmmm}) that $H_{cov}(X,x_0)=h(x_0)$, using Lemma \ref{moment de
psi}. A rigourous proof of this statement is given in Section \ref{covproof}.
Moreover,   the  continuity  of   $h$  around   $y_0+t\alpha$,  for   $y_0  \in
<\alpha>^{\perp}$  and  $t\in\R$,  applied  in  (\ref{covxraymmm})  and  Lemma
\ref{GdansL2} will imply that
$H_{cov}({\mathcal P}_{\alpha}X,y_0)=m(\alpha, y_0)+1/2$,
where
\begin{equation}\label{minimum}
m(\alpha,  y_0):=\inf_{t\in\R}h(y_0+t\alpha).
\end{equation}
These observations yield to the next theorem. A detailed proof is
given  in Sections \ref{covproof} and \ref{fddproof}. We denote {\sl meas} the Lebesgue's measure.
\begin{thm} \label{thlassh}
Let $d\ge 2$.  Let $h$ be a $\beta$-Lipschitz  function  on  $\R^{d}$ such
that $0<  h < 1/2$. Let $X$ be the multifractional microball model with index $h$
and let  ${\mathcal P}_{\alpha}X$  be its $X$-ray  transform in  the direction
$\alpha \in S^{d-1}$.
\begin{itemize}
\item At any point
$x_0\in \R^d$,  the cov-lass index of $X$ is equal to $h(x_0)$ and the fdd-lass index of $X$ is greater than $\min(1,\beta+2h(x_0))$. Moreover the
covariance  of   $\lambda^{-h(x_0)}\Delta_{x_0}X(\lambda.)$  converges,  up   to  a
multiplicative constant, to the
covariance of a fractional Brownian motion of index $h(x_0)$.
\item  At any  point $y_0\in  <\alpha>^{\perp}$,  the cov-lass  index and  the
  fdd-lass index of ${\mathcal P}_{\alpha}X$  are equal to $m(\alpha,y_0)+1/2$.  Moreover,
  when  meas$(\{t\in \R^*;~~h(y_0+t\alpha)=m(\alpha, y_0)\})>0$, the
covariance and the finite dimensional distributions of
$\lambda^{-m(\alpha,    y_0)+1/2}\left(\Delta_{y_0}({\mathcal    P}_{\alpha}X)
  (\lambda.)-\mathbf{E}\left(\Delta_{y_0}({\mathcal               P}_{\alpha}X)
  (\lambda.)\right)\right)$
converge, up to a
multiplicative constant, to the corresponding ones  of a fractional Brownian motion of index $m(\alpha, y_0)+1/2$.
\end{itemize}

\end{thm}
\begin{rem} The first point of Theorem \ref{thlassh} is still true in the one-dimensional case ($d=1$).
\end{rem}

 We remark that the cov-lass index of $X$ at point $x_0\in\R^d$  is
 equal to $h(x_0)$, as for the
multifractional Brownian motion with index $h$, and this justifies
the name  multifractional. Also
  let us just point out that these results apply to the fractional case, where
 $h$ is constant. In this case we have shown that the
fdd-lass index for the microball is  equal to $+\infty$. This is a consequence
 of the stationarity of the model. The first
order moment of the increments is then equal to $0$ and this is no
lomger the case when $h$ is not constant.
 From  a certain point of view, for the microball model, the fdd-lass index deals
 with the regularity of the first order moment
and the cov-lass index with those of the second order moment.
 Finally, compared to the fractional case,
 we still have an additive factor of $1/2$ for the lass indices of the
 windowed X-ray  transforms, but here only  the infimum of  $h$ along straight
 lines can be recovered. \\

\subsubsection{The singular case:}

  let us now consider $h$  to be  an even,  $\beta$-Lipschitz function  on the
sphere, defined on $\R^d\smallsetminus\{0\}$ by homogeneity.  Of
course, this case  is only relevant when  $d\ge 2$ because otherwise
it turns to be the fractional case, with $h$ constant. The
singularity of $h$ at point 0 makes this point a very special one.
The balls are thrown from 0 and their
  numbers and sizes only depend on the direction along which they are thrown.
  To distinguish the microball model associated with such a singular $h$ from the multifractional one, we will call it the star microball model. 
Let us remark that the $\beta$-Lipschitz assumption on $h$ means that there exists $C>0$ such that for all $x_0\in\R^d\smallsetminus\{0\}$ and $\xi\in \R^d$, when $\m\xi\m\le 1$
\begin{equation}\label{betalip}
 \left|h(x_0+\xi)-h(x_0)\right|\le C\m x_0\m^{-\beta}\m\xi\m^{\beta}.
\end{equation}

 Then,  the lass properties of the star microball model, respectively its X-ray transform in the direction $\alpha\in S^{d-1}$, at any point $x_0\in \R^{d}\smallsetminus\{0\}$,  respectively $y_0\in<\alpha>^{\perp}\smallsetminus\{0\}$, are the same as for the multifractional microball model, respectively its X-ray transform in the direction $\alpha$, given in Theorem \ref{thlassh}.  Thus, the next theorem will only deal with the lass properties around $0$. Let us remark that, in that case, by homogeneity of $h$, the exponent in
(\ref{covmmm}) is equal to $\lambda^{h(\xi)}$ for all $\xi\neq 0$.
Then, denoting by $m\in (0, 1/2)$ the minimum  of $h$ on $S^{d-1}$,  by Lemma \ref{moment de psi}, we will get
    $H_{cov}(X,0)=m$.
         Moreover, in (\ref{covxraymmm}),  since  $h$  is continuous around   $t\alpha$ and  $h(t\alpha)=h(\alpha)$,  for    $t\in\R^*$, using Lemma \ref{GdansL2}, we prove that $H_{cov}({\mathcal P}_{\alpha}X,0)=h(\alpha)+1/2$.\\
Let us state the different lass properties of the star microball model.
\begin{thm} \label{thlasshstar}  Let $h$ be an even $\beta$-Lipschitz, non constant,  function on  $S^{d-1}$ such
that $0<  h < 1/2$.
Let $X$ be the star microball model with index $h$
and let  ${\mathcal P}_{\alpha}X$  be its $X$-ray  transform in  the direction
$\alpha \in S^{d-1}$. We consider lass properties at point $0$.
\begin{itemize}
\item The cov-lass index of $X$ is equal to $m$ and the fdd-lass index of $X$ is equal to $2m$. Moreover, when  meas$(\{h=m\})>0$, the
covariance  of   $\lambda^{-m}\Delta_{0}X(\lambda.)$  converges to $\gamma_m$ with
$$\gamma_m(x,x')=\int_{\R^d\times \R^+}\mathbf{1}_{h(\xi)=m}\psi(
x,\xi,r)\psi( x',\xi,r)r^{-d-1+2m}d\xi dr~,$$ for $ x, x'\in \R^d,$
while        the       finite        dimensional
distributions       of\\
$\lambda^{-2m}\left(\Delta_{0}X(\lambda.)-\mathbf{E}\left(\Delta_{0}X(\lambda.)\right)\right)$
converge to the deterministic field $Z_m$, with
$$Z_m=\left\{-\m x\m ^{2m}\int _{\R^d\times\R^+} {\bf 1}_{h(\xi)=m}
  \psi(\frac{x}{\m x \m},\xi,r) r^{-d-1+2m}~ d\xi dr~; ~x\in \R^d~\right\},$$
  where $\psi$ is given by (\ref{psi}).
\item  The cov-lass  index and  the
  fdd-lass index of ${\mathcal P}_{\alpha}X$  are equal to $h(\alpha)+1/2$.  Moreover,
   the
covariance and the finite dimensional distributions of $\lambda^{-h(\alpha)+1/2}\left(\Delta_{0}({\mathcal P}_{\alpha}X) (\lambda.)-\mathbf{E}\left(\Delta_{0}({\mathcal P}_{\alpha}X) (\lambda.)\right)\right)$ converge, up to a
multiplicative constant, to the corresponding ones  of a fractional Brownian motion of index $h(\alpha)+1/2$.
\end{itemize}

\end{thm}

The proof is given in Sections \ref{covproof} and \ref{fddproof}. Let us remark that  for the star
microball model,  there  exists both a cov-lass index and a fdd-lass
index  and  that  the  later   equals  the  double  of the former.
This   multiplicative  factor   is  typical  for   the Poisson
structure  proved by the following exercise.

\begin{rem}
Let $(X_n)$ be a sequence of Poisson random variables. Suppose there exists
some $H>0$ and $v>0$ such that
$\mbox{Var}(n^H(X_n-E(X_n)))$ tends  to $v$ when $n$ tends  to $+\infty$. Then
$n^{2H}(X_n-E(X_n))$ tends in distribution to $-v$.
\end{rem}

Moreover when  $\{h=m\}$ has positive measure, the tangent field at
0 is deterministic and not zero, hence does not have stationary
increments. This is worth to be noticed and linked to a result of
Falconer \cite{Falc}, which states that at almost all points the
tangent field -if it exists-  must have stationary increments. Hence
the point $0$ appears as an `exceptional point'  (see \cite{CL2} for
other examples of exceptional points).

Finally, let us point out that the tangent field of the X-ray
transform, when it exists, is Gaussian, even a fractional Brownian
motion, whereas the tangent field of the star microball model was
deterministic. This  justifies, from a mathematical point of view,
modeling radiographic images by fBm, even when the media  under
study is far from being of this type (see \cite{Harba} for an
experimental study).

\subsection{Comparison with homogenization results} \label{sshom}

There  are different  ways  to  consider     self-similarity at
small scales, depending on which part of the signal the scaling
acts.  Instead of performing a  scaling   on  the   increments  lag,
as   it  is  performed   in  Sections \ref{ssfmm} and \ref{ssmfmm},
we act on the radius of the balls  as follows. Suppose we  zoom and
consider the balls $B(\xi,r/\varepsilon)$ instead of  the balls
$B(\xi,r)$, where  the $(\xi,r)$ are  randomly chosen by the Poisson
random measure $N_h$, and we let $\varepsilon$ tend to 0. Denoting
by $X^{\varepsilon}$ the associated field
$$X^{\varepsilon}(x)=\int_{\R^d\times
  \R^+}\mathbf{1}_{B(\xi,r/\varepsilon )}(x)N_h(d\xi,dr)~,~ x\in \R^d~,$$
we look for  normalization terms $n(x_0)$ such that
$\epsilon^{n(x_0)}(\Delta_{x_0}X^{\varepsilon}-\mathbf{E}(\Delta_{x_0}X^{\varepsilon})$
converge  in distribution to  a non degenerate field.  Note that the
field $X^{\varepsilon}$ can also be considered  as a microball model
(see Definition \ref{def}) associated with a Poisson measure with
intensity
$$\nu^{\varepsilon}_h(d\xi,dr)=\varepsilon^{-d+2h(\xi)}r^{-d-1+2h(\xi)}\mathbf{1}_{(0,\varepsilon^{-1})}(r)d\xi dr~.$$
Actually this procedure is nothing but homogenization and is
 close to the
thermodynamical  limit investigated in  \cite{Mand1} or  the scaling
limit in \cite{KAJ}. Similar computations as for  the previous
theorems yield\\
$\bullet$ if $0<m:=\min h<1/2$ and the set $\{\xi;h(\xi)=m\}$ has
positive measure, then the normalization term $n(x_0)$ is equal to
$d/2-m$ for all $x_0\in \R^d$\\
$\bullet$ moreover if $h$ is constant equal  to $m$, then the limit field is a
fractional Brownian motion with index $m$.


\section{Conclusion}

We propose to model -from a microscopic point of view- the mass intensity of a porous media or the number of
connected custumers in  a network with a non-Gaussian  field, which presents a
macroscopic (multi)fractional behavior. The
rich structure of Poisson point  processes allows us to reach this
goal and also to perform
explicit computations  as  in  the
Gaussian case.  In order to  keep the model as  intuitive as
possible, we did not try to produce more general fields. The
Poisson structure can  obviously be exploited  further on by
considering more  general integrators with  respect  to  the
Poisson  measure.  Replacing  the  indicator  function
$\mathbf{1}_{B(\xi,r)}$  in Definition \ref{def}  by a  more
general  one will, for instance, give the possibility to model
granular media with non spherical  grains. Another model
for porous media  can also be built up from a collection of
random spheres  which correspond no  more  to grains,
but  to pores  or bubbles. By this way, one will  get a
$\{0,1\}$-valued field and leave the linear context.

\section{Proofs of the cov-lass properties}\label{covproof}
Let  $h$ be a function defined  on $\R^d$ or $\R^d\smallsetminus\{0\}$, with
$0<h <1/2$.\\
In this section we will  give rigorous proofs for the cov-lass properties of the microball models and their
windowed X-ray transforms.
\subsection{Cov-lass properties of the microball model}\label{covlassmicroball}
 Let $x_0\in \R^d$. For $H\in (0,1)$ and $x, x'\in \R^d$, recall from (\ref{covmmm}) that
$$\Gamma^H_{\lambda}(x_0,x,x'):=\mbox{Cov}\left(
\frac{\Delta_{x_0}X(\lambda x)}{\lambda^{H}},\frac{\Delta_{x_0}X(\lambda x')}{\lambda^{H}}\right)$$
$$=
\int_{\R^d\times
(0,\lambda^{-1})}\lambda^{-2(H-h(x_0+\lambda\xi))}\psi( x,\xi,r)\psi(
x',\xi,r)r^{-d-1+2h(x_0+\lambda\xi)}d\xi dr.$$
\\

a. {\it The smooth case:} we assume that there exists $C>0$ such that, for $\m\xi\m<1$
\begin{equation}\label{betalipmmm}
\m h(x_0+\xi)-h(x_0)\m\le C\m\xi\m^{\beta}.
\end{equation}
 We  first establish that
\begin{equation} \label{conv de gamma}
\Gamma^{h(x_0)}_{\lambda}(x_0,x,x')\stackrel{\lambda\rightarrow 0^+}{\longrightarrow}
\Gamma^{h(x_0)}(x,x')~,
\end{equation}
where
$$\Gamma^{h(x_0)}(x,x')=\int_{\R^d\times \R^+}\psi(
x,\xi,r)\psi( x',\xi,r)r^{-d-1+2h(x_0)}d\xi dr ~,$$
and  then prove  that  $\Gamma^{h(x_0)}$  is the  covariance function of a  fractional
Brownian motion with Hurst index $h(x_0)$.

Since $h(x_0)<1/2$,  Lemma \ref{moment de  psi} indicates that  the function
$\psi(x,.)$ belongs to $L^2(\R^d\times \R^+, r^{-d-1+2h(x_0)}d\xi dr)$.
By Cauchy-Schwarz inequality, in order to prove (\ref{conv de gamma}), it is enough to prove that the difference
$I(\lambda)$, given by
\begin{equation}\label{difference1c}
I(\lambda)=\int_{\R^d\times (0,\lambda^{-1})}\psi(
x,\xi,r)^2\left|(\lambda r)^{-2(h(x_0)-h(x_0+\lambda\xi))}-1\right| r^{-d-1+2h(x_0)}d\xi dr~,
\end{equation}
tends to $0$ as $\lambda$ tends to $0^+$.
Let us remark that, by (\ref{betalipmmm}),  for some positive $p$ and $q\in(0,1)$ to be fixed later, when $\lambda^{p+1}<\lambda r<1$ and $\lambda\m\xi\m<\lambda^{1-q}<1$,
 one can find $C>0$ such that
 \begin{equation}\label{majorationpremierterme}
\left|(\lambda r)^{-2(h(x_0)-h(x_0+\lambda\xi))}-1\right|\le
C\left(\lambda\m\xi\m\right)^{\beta}\m\ln{\lambda}\m.
\end{equation}
Without loss of generality, we can assume $\beta<1-2h(x_0)$. In this case
we will prove that, for  $\lambda$ small enough compared to $x$,
\begin{equation}\label{majorationdifference}
I(\lambda)~\le C\lambda^{\beta}\m\ln{\lambda}\m,
\end{equation}
with $C>0$.
We split the integral into
$$\int_{\R^d\times  (0,\lambda^{-1})}=\int_{B(0,\lambda^{-q})\times
  (\lambda^{p},\lambda^{-1})}
+ \int_{\R^d\times (0,\lambda^{p})}
+\int_{B(0,\lambda^{-q})^c\times (0,\lambda^{-1})}.$$
By (\ref{majorationpremierterme}), the same kind of arguments as in the proof of Lemma \ref{moment de psi} yield to
$$\left|\int_{B(0,\lambda^{-q})\times
  (\lambda^{p},\lambda^{-1})}\right|\le C\lambda^{\beta}\m x\m^{\beta+2h(x_0)}\m\ln{\lambda}\m.$$
Moreover, for $\xi\in B(0,\lambda^{-q})^c$ or $r\in (0,\lambda^{p})$, we use the following inequality
$$\left|(\lambda    r)^{-2(h(x_0)-h(x_0+\lambda\xi))}-1\right|\le    2(\lambda
r)^{2(m-h(x_0))},$$
which holds since $h(x_0+\lambda
\xi)\ge m$ and $\lambda r\le 1$. Then, on one hand,
\begin{eqnarray*}
\left|\int_{\R^d\times (0,\lambda^{p})}\right|&\le &2\lambda^{2(m-h(x_0))}\int_{\R^d\times (0,\lambda^{p})}\mathbf{1}_{\m \xi\m <r}r^{-d-1+2m}dr\\
&\le &C\lambda^{2(m-h(x_0))+2mp}.
\end{eqnarray*}
On the other hand, for $\lambda\le \left(4(1+\m x\m)\right)^{-1/q}$,
\begin{eqnarray*}
\left|\int_{B(0,\lambda^{-q})^c\times (0,\lambda^{-1})}\right|&\le & C\m x\m\lambda^{2(m-h(x_0))}\int_{\m\xi\m>\lambda^{-q}}
\m\xi\m^{-d-1+2m}d\xi\\
&\le& C\m x\m \lambda^{2(m-h(x_0))+q(1-2m)}.
\end{eqnarray*}
Finally,  it suffices to choose $p$ and $q$ such that
$$p>\frac{\beta+2(h(x_0)-m)}{2m}\,\, \mbox{ and } \,\,\frac{\beta+2(h(x_0)-m)}{1-2m}<q<1.$$
This concludes for the proof  of the convergence  of $\Gamma^{h(x_0)}_{\lambda}(x_0,.)$
to $\Gamma^{h(x_0)}$.  It remains to  show that $\Gamma^{h(x_0)}$  is the
covariance of a -up to a constant- fractional Brownian motion with Hurst index
$h(x_0)$.
A straightforward computation gives
\begin{eqnarray*}
\int_{\R^d\times \R^+} && \left(
  \psi(x,\xi,r)-\psi(x',\xi,r)\right)^2 ~ r^{-d-1+2h(x_0)}~d\xi dr \\
&=& \int_{\R^d\times \R^+}    \psi(x-x',\xi,r)^2 ~ r^{-d-1+2h(x_0)}~d\xi dr
 ~.
\end{eqnarray*}
This allows to write $\Gamma^{h(x_0)}$ as
$$\Gamma^{h(x_0)}(x,x')=\frac{1}{2}\left(v(x)+v(x')-v(x-x')\right)$$
where, $v(x)=\Gamma^{h(x_0)}(x,x)=c|x|^{2h(x_0)}.$
\\
This proves the covariance part of the first point in Theorem \ref{thlassh}.

b.  {\it The  singular case:}  we now  assume  that $h$  is given  by an  even
$\beta$-Lipschitz function on the sphere, defined on $\R^d\setminus\{0\}$ by homogeneity,
and deal with the case $x_0=0$.\\
b-i) Assume that $\{\xi \in \R^d\setminus\{0\};h(\xi)=m\}$ has positive measure.\\
Thus, using Lemma \ref{moment de psi}, Lebesgue's Theorem gives that $\Gamma^{m}_{\lambda}(0,x,x')$ tends to
\begin{equation}\label{limfract}
\int_{\R^d\times \R^+}\mathbf{1}_{h(\xi)=m}\psi(
x,\xi,r)\psi( x',\xi,r)r^{-d-1+2m}d\xi dr~
:=\gamma^{m}(x,x'),
\end{equation}
which does not vanish by assumption. \\
\\
b-ii) Assume now that  $\{\xi \in \R^{d}\setminus\{0\}; h(\xi)=m\}$ has
measure $0$.\\ Since $\gamma^{m}(x,x')=0$ by (\ref{limfract}), the covariance
$\Gamma^H_{\lambda}(0,x,x')$  tends to $0$ for $H\le m$.
On the other hand, for $H=m+2\epsilon$ with $\epsilon>0$ and for
$\lambda$ in $(0,1)$,
$$\Gamma^{H}_{\lambda}(x,x)\ge \lambda^{-\epsilon}
\int_{\R^d\times (0,1)}\mathbf{1}_{h(\xi)<m+\epsilon}\psi( x,\xi,r)^2~r^{-d-1+2h(\xi)}d\xi dr~,$$
where $\{\xi \in \R^d\setminus\{0\}; h(\xi)<m+\epsilon\}$ has positive measure.
So the above quantity tends to infinity when $\lambda$ tends to
$0^+$. Hence the exponent $m$ is proved to be the  cov-lass index for
$X$ at $x_0=0$.\\ \\
The   cov-lass   properties  asserted   in   the   first   part  of   Theorem
\ref{thlasshstar} are now proved.
\medskip

\subsection{Cov-lass properties of the windowed X-ray transforms}
Let $\alpha\in S^{d-1}$ and $y_0\in <\alpha>^{\perp}$.
Let us denote by  $\Gamma_{\lambda}^{H}(y_0,.)$ the covariance function of the
field $\lambda^{-H}\Delta_{y_0}{\mathcal
P}_{\alpha}X(\lambda .).$
Then for $y$ and $y'$ in $<\alpha>^{\perp}$, by (\ref{covxraymmm})
\begin{eqnarray*}
\Gamma_{\lambda}^{H}(y_0,y,y')=\lambda^{-2H+1}&&
\int_{<\alpha>^{\perp}\times \R\times (0,\lambda^{-1})}
G_{\rho(t+\lambda .)}(y,\gamma,r)G_{\rho(t+\lambda .)}(y',\gamma,r)\\
&& \times ~(\lambda r)^{2h(\lambda \gamma+t\alpha+y_0)}r^{-d-1}~d\gamma dt dr.
\end{eqnarray*}
Let us write $G_{\rho(t+\lambda .)}$ as the integral given by (\ref{grho})
\begin{eqnarray*}
\Gamma_{\lambda}^{H}(y_0,y,y')=\lambda^{-2H+1}&&
\int_{<\alpha>^{\perp}\times \R\times (0,\lambda^{-1})\times\R}
\psi(y,\gamma-p\alpha,r)\rho(t+\lambda p)(y,\gamma,r)\\
&& \times ~G_{\rho(t+\lambda .)}(y',\gamma,r) (\lambda r)^{2h(\lambda \gamma+t\alpha+y_0)}r^{-d-1}~d\gamma dt dr dp.
\end{eqnarray*}
Two  changes of variables allow to write  $\Gamma_{\lambda}^{H}(y_0,y,y')$ as
$$\lambda^{-2H+1}
\int_{ \R\times (0,\lambda^{-1})\times\R^d}
\psi(y,\xi,r)\rho(t)G_{\rho(t+\lambda .)}(y',\xi,r)
 \times ~(\lambda r)^{2h(\lambda \xi+t\alpha+y_0)}r^{-d-1}~ d\xi dr dt.
$$

a.  {\it The smooth case:} we assume that $h$ is $\beta$-Lipschitz.
By the same  arguments as in the proof  of (\ref{majorationdifference}), using
the fact that $G_{\rho}\le 2\n\rho\n_{\infty}\m G\m$,
for $\lambda$ small enough compared to $y$ and $y'$,
$$\int_{ \R^d\times (0,\lambda^{-1})}
\psi(y,\xi,r)G_{\rho(t+\lambda .)}(y',\xi,r)
 \left|(\lambda r)^{2h(\lambda \xi+t\alpha+y_0)-2h(t\alpha+y_0)}-1\right|$$
\begin{equation}\label{cvgce}
\times~ r^{-d-1+2h(t\alpha+y_0)}~d\xi dr \le C\lambda^{\beta}\m\ln(\lambda)\m,
\end{equation}
with $C>0$ and $\beta<1-2h(t\alpha+y_0)$.
Thus we write $(\lambda    r)^{2h(\lambda    \xi+t\alpha+y_0)}$ as
$$(\lambda   r)^{2h(\lambda
  \xi+t\alpha+y_0)-2h(t\alpha+y_0)}\times
(\lambda       r)^{2h(t\alpha+y_0)-2m(\alpha,y_0)}       \times       (\lambda
  r)^{2m(\alpha,y_0)},$$
  where  $m(\alpha,y_0)$ is given by (\ref{minimum}).
For almost every $t\in\R$, the second factor tends
to $\mathbf{1}_{h(y_0+t\alpha)=m(\alpha,y_0)}$.
Since $m(\alpha,y_0)<1/2$, let us recall that by    Lemma   \ref{GdansL2} the function $G(y,.)$  belongs to $L^2(<\alpha>^{\perp}\times
\R^+,r^{-d-1+2m(\alpha, y_0)}d\gamma    dr)$. We use Lebesgue's Theorem and (\ref{cvgce}),  to obtain the following asymptotics:
$\Gamma_{\lambda}^{m(\alpha, y_0)+1/2}(y_0,y,y')$ tends to
$$\left( \int_{\R}\mathbf{1}_{h(y_0+t\alpha)=m(\alpha,y_0)}
\rho^2(t)dt \right)\times \Gamma^{m(\alpha, y_0)+1/2}(y,y'),$$
where
\begin{equation} \label{Gamma}
\Gamma^{H}(y,y'):=
\int_{<\alpha>^{\perp}\times\R^+} G(y,\gamma,r)G( y',\gamma,r)
   r^{-d-2+H}d\gamma dr.
\end{equation}
 The
identification of $\Gamma^{H}$ as the covariance of a fractional Brownian
motion defined on $<\alpha>^{\perp}$ with Hurst index
$H=m(\alpha,y_0)+1/2$ is straightforward following the same
arguments as in the part a. of  \ref{covlassmicroball}.

\medskip

If $\{t\in  \R;h(y_0+t\alpha)=m(\alpha, y_0)\}$ has positive  measure, then we
have finished  with the proof.  Otherwise, we proceed  in the same way  as in
part  b-ii) of  \ref{covlassmicroball}. The  cov-lass property  of  the second
part of Theorem \ref{thlassh} is established.\\

b.  {\it The singular case:} we assume that
 $h$  is  $\beta$-Lipschitz  on  $S^{d-1}$  with
$0<\beta\le 1$ and look at the cov-lass property around 0. By homogeneity,
the $\beta$-Lipschitz condition is replaced by
$$\left|h(t\alpha+x)-h(\alpha)\right|\le Ct^{-\beta}\m x\m^{\beta}.$$
Hence the upper bound of (\ref{cvgce}) is now given by
$$
 C\left(\m t\m^{-(1-2m)}+\m t\m^{-(1-2h(\alpha))}\right)\lambda^{\beta}\m\ln(\lambda)\m.$$
Lebesgue's theorem still applies to get  the result enounced in the second part 
 of Theorem \ref{thlasshstar}.

\section{Fdd-lass properties}\label{fddproof}
We now prove   the finite dimensional distribution lass
properties.
\subsection{Fdd-lass properties of the microball model}
 Let us denote $\widetilde{X}=X-\mathbf{E}\left(X\right)$ the centered version of $X$. For notational sake of
simplicity, we will only consider the limit in distribution of
$\lambda^{-H} \Delta_{x_0}\widetilde{X}(\lambda x)$
 for a fixed $x$ in
$\R^d$ instead of a  random
vector   $(\lambda^{-H} \Delta_{x_0}\widetilde{X}(\lambda x_j))_{1\le  j \le  n}$. The general case follows along the same lines.

For $H>0$, $x\in\R^d$ and $t\in\R$, let
$$\mathbf{E}\exp \left(it \frac{\Delta_{x_0}\widetilde{X}(\lambda x)}
{\lambda^{H}} \right)
=\exp \Phi (H,\lambda,x_0,x,t)$$
where $\Phi (H,\lambda,x_0,x,t)$ is given by
$$
\int_{\R^d\times\R^+} \left( e^{it\lambda^{-H} \psi(\lambda x,\xi-x_0,r)}
-1-i t\lambda^{-H} \psi(\lambda x,\xi-x_0,r) \right)
\,d\nu_h(\xi,r)~.$$
A change of variable yields
\begin{eqnarray*}
\Phi(H,\lambda,x_0, x, t)
&=&\int_{\R^d\times (0,\lambda^{-1})} \left( e^{it\lambda^{-H} \psi(x,\xi,r)}
-1 -i t\lambda^{-H} \psi(x,\xi,r) \right)
  \\
&\times& ~\lambda^{2h(x_0+\lambda\xi)}~r^{-d-1+2h(x_0+\lambda\xi)}~d\xi dr.
\end{eqnarray*}
Lemma \ref{moment de psi} allows us to split the integral into
$\Phi=\Phi_1+(\Phi-\Phi_1)$, where $\Phi_1 (H,\lambda,x_0, x,t)$ is equal to
$$
\int_{\R^d\times(0,\lambda^{-1})} \left( e^{it\lambda^{-H} \psi(x,\xi,r)}
-1 \right) \lambda^{2h(x_0+\lambda\xi)} ~r^{-d-1+2h(x_0+\lambda\xi)}
~d\xi dr.$$
\noindent
Then,
$$|\Phi_1 (H,\lambda,x_0, x, t)| \leq \lambda^{2m} \int_{\R^d\times\R^+}
\left| e^{it\lambda^{-H} \psi(x,\xi,r)}-1 \right|  r^{-d-1+2m}~ d\xi dr~.$$
We notice that
$$\left| e^{it\lambda^{-H} \psi(x,.)}-1 \right|
\leq  2 ~\ind_{t \psi(x,.) \ne 0}
\leq  2 ~ |\psi(x,.)|$$
and recall that $\psi(x,.)$ belongs to $L^1(\R^d\times \R^+,r^{-d-1+2m}~dr
d\xi)$ by Lemma \ref{moment de psi} so that
$$ \lim _{\lambda \rightarrow 0^+}\Phi_1 (H,\lambda,x_0, x,t)=0.$$
The second term $\Phi_2:=\Phi-\Phi_1$ is given by
\begin{equation} \label{phi2}
-it~ \int _{\R^d\times(0,\lambda^{-1})}
\lambda^{-H+2h(x_0+\lambda\xi)}  \psi(x,\xi,r)
 ~ r^{-d-1+2h(x_0+\lambda\xi)}~d\xi dr~.
\end{equation}
We will now analyse on this expression.

First, let us note that in the case where $h$ is constant,
 then the increments of $X$ are equal to the increments of $\widetilde{X}$ by stationarity and
 hence $\Phi_2=0$. That concludes the proof for the fractional
microball model, i.e.  the first part of Theorem \ref{thlassm}.

Now we
deal with the multifractional model (case $h$ non constant).\\ \\
a. {\it The smooth case:} We assume that (\ref{betalipmmm}) holds and
 we prove that the critical fdd-lass
index of $X$ at $x_0$ is  $\ge \min(1,\beta+2h(x_0))$. Without loss of generality, we can assume that $\beta+2h(x_0)\le 1$.
For
$H<\beta+2h(x_0)$, we will establish that $\Phi_2 (H,\lambda,x_0, x, t)$ tends to $0$ when $\lambda\rightarrow
0^+$.
Recall that by equation (\ref{phi2})  $\Phi_2(H,\lambda,x_0, x, t)$ is equal to
$$-it~\lambda^{-H+2h(x_0)}
\int _{\R^d\times(0,\lambda^{-1})}
(\lambda r)^{2(h(x_0+\lambda\xi)-h(x_0))}  \psi(x,\xi,r)
 ~ r^{-d-1+2h(x_0)}~d\xi dr~.
$$
Since $\int _{\R^d\times(0,\lambda^{-1})}
 \psi(x,\xi,r) ~ r^{-d-1+2h(x_0)}~d\xi dr~$ vanishes, we get
$$
\Phi_2 (H,\lambda,x_0, x, t)
=-it~\lambda^{-H+2h(x_0)} $$
$$\times \int _{\R^d\times(0,\lambda^{-1})}
\left( (\lambda r)^{2(h(x_0+\lambda\xi)-h(x_0))}-1 \right)  \psi(x,\xi,r)
 ~ r^{-d-1+2h(x_0)}~d\xi dr.$$

By the $\beta$-Lipschitz assumption on $h$,
for $\lambda$ small enough compared to $x$, using (\ref{majorationdifference}) and since $\m\psi\m=\psi^2$, we get
$$\int _{\R^d\times(0,\lambda^{-1})}
\left| (\lambda r)^{2h(x_0+\lambda\xi)-2h(x_0)}-1 \right|  |\psi(x,\xi,r)|
 ~ r^{-d-1+2h(x_0)}~d\xi dr \le C(x)\lambda^{\beta}.$$
This implies the lower bound for the fdd-lass index. \\
\\
b. {\it The singular case:} We now assume that $h$ is  $\beta$-Lipschitz on the sphere and only deal with the case $x_0=0$.\\
b-i) Assume that $\{\xi;h(\xi)=m\}$ has
 positive measure.
 For $H=2m$, using Lemma \ref{moment de psi} again, we get
\begin{equation}\label{limmultifdd} \lim _{\lambda \rightarrow 0^+}\Phi_2 (2m,\lambda,0, x,t)= -it~\int _{\R^d\times\R^+} {\bf 1}_{h(\xi)=m}
  \psi(x,\xi,r) r^{-d-1+2m}~ d\xi dr~.
  \end{equation}
 Hence the finite dimensional distributions of $\lambda^{-2m}\Delta_0X(\lambda .)$
converge to the finite dimensional distributions
of the  deterministic field $Z^m$ given by
$$Z^m(x)=-\int _{\R^d\times\R^+} {\bf 1}_{h(\xi)=m}
  \psi(x,\xi,r) r^{-d-1+2m}~ d\xi dr~, ~x\in \R^d~.$$
It remains to show that $Z^m$ is not zero. This follows from the next lemma,
  since $\{\xi;h(\xi)\ne m\}$ contains a ball by continuity of $h$.

\begin{lem}\label{zeroZm}
Let $m\in (0,1/2)$.
For all Borel sets $E\subset \R^d$ with positive measure, if
$$\int_{\R^d\times \R^+}\mathbf{1}_E(\xi)\psi(.,\xi,r)r^{-d-1+2m}d\xi dr =0$$
 in $\R^d$,
then the set $E^c$ does not contain any open ball.
\end{lem}
\begin{proof}

Note that for all $x\in\R^d$,
$$\int_{\R^+}\psi(x,\xi,r)r^{-d-1+2m}d\xi dr= (2m-d)^{-1} \left(
|x-\xi|^{-d+2m} - |\xi|^{-d+2m} \right)$$
and hence, for all Borel sets $E\subset\R^d$,
\begin{eqnarray*}
I_E^m(x):&=&\int_{\R^d\times \R^+}\mathbf{1}_E(\xi)\psi(x,\xi,r)r^{-d-1+2m}d\xi
dr\\
&=& \int_{\R^d}\mathbf{1}_E(\xi)(2m-d)^{-1} \left(
|x-\xi|^{-d+2m} - |\xi|^{-d+2m} \right) d\xi ~.
\end{eqnarray*}
Let us suppose that we can find an open nonempty ball $B\subset E^c$. Then
 $I_E^m$ is smooth on $B$. The Laplacian of
$I_E^m$ can easily be computed
$$\Delta I_E^m(x)=-2(1-m)\int_{\R^d}\mathbf{1}_E(\xi)
|x-\xi|^{-d+2m-2} d\xi ~,$$
and proved to be negative on $B$. Thus, $I_E^m$ does not vanish on $B$. This completes the proof of Lemma \ref{zeroZm}.
\end{proof}

\medskip

\noindent b-ii) Assume that $\{\xi;h(\xi)=m\}$
is of measure  $0$. we will establish that  the fdd-lass index of $X$  at 0 is
still equal to $2m$.   For  $H\leq 2m$, from (\ref{limmultifdd}), $\Phi_2
(H,\lambda,0,x,t)$ tends to 0.
On the other hand, for  $H\in (2m,1)$ and $H<2 \max h$, we
have to prove that there exists at least one $x\in\R^d$ such
that  $\Phi_2 (H,\lambda,0,x,t)$ does not tend to 0.
First, let us remark that, for all $x\in\R^d$, $\Phi_2 (H,\lambda,0,x, t)$ may be written as
$$
-it~ \left(\int _{\R^d\times\R^+}
\lambda^{-H+2h(\xi)}  \psi(x,\xi,r)
 ~ r^{-d-1+2h(\xi)}~d\xi dr~+{\mathcal O}\left(\lambda^{1-H}\right)\right).$$
Thus it is sufficient to consider
\begin{eqnarray*}
\tilde{\Phi}_{2}(H,\lambda,x):
&=& \int _{\R^d\times\R^+}
\lambda^{-H+2h(\xi)}  \psi(x,\xi,r)
 ~ r^{-d-1+2h(\xi)}~d\xi dr\\
&=& \int _{\R^d} (2h(\xi)-d)^{-1}
 \lambda^{-H+2h(\xi)}
f(x,\xi)  ~d\xi~,
\end{eqnarray*}
with
$$f(x,\xi)=|x-\xi|^{-d+2h(\xi)} - |\xi|^{-d+2h(\xi)}.$$
The function $f(.,\xi)$ is smooth on $\R^d\smallsetminus\{\xi\}$ and its Laplacian given by
$$\Delta f(x,\xi)=(2h(\xi)-d)(h(\xi)-1)|x-\xi|^{-d+2h(\xi)-2},$$
is negative. Then, approximating the Laplacian by the second order increments, one can find $C\in (0,\frac{1}{2})$ such that, whenever $\m x-\xi\m \neq 0$ and
 $\delta\le C\m x-\xi\m$,
\begin{equation}\label{Laplacien}
\sum_{1\le j \le
d} f(x+\delta e_j,\xi)+f(x-\delta
e_j,\xi)-2f(x,\xi)\le \frac{\delta^2}{2}\Delta f(x,\xi).
\end{equation}
Let us take $H$ such that $2m<H<2~\max h$ and $H<1$
and note that the sets $\{\xi;H\le 2h(\xi)\}$ and
$\{\xi;H>2h(\xi)\}$ have positive measure. Then, by continuity of $h$,
there exists a nonempty open ball $B\subset \{\xi;H>2h(\xi)\}$. For every $x\in\R^d$ we introduce
$$\tilde{\Phi}_{\lambda}(x):=\int _{\R^d} {\bf 1}_{2h(\xi)\le H}(2h(\xi)-d)^{-1}
 \lambda^{-H+2h(\xi)}
f(x,\xi)  ~d\xi,$$
such that
$$\tilde{\Phi}_{2}(H,\lambda,x)=\tilde{\Phi}_{\lambda}(x)+\left(\tilde{\Phi}_{2}(H,\lambda,x)-\tilde{\Phi}_{\lambda}(x)\right).$$
By Lebesgue's Theorem, the second term tends to 0 with $\lambda$.
Suppose that $\tilde{\Phi}_{\lambda}(x)$ tends to 0 with $\lambda$
for every $x$ in $\R^d$. Then for all $x\in \R^d$ and all $\delta
\in \R$,
\begin{equation} \label{delta2phi}
\Delta ^{(2)}_{\delta}\tilde{\Phi}_{\lambda}(x):=\sum_{1\le j \le
d} \left( \tilde{\Phi}_{\lambda}(x+\delta e_j)+\tilde{\Phi}_{\lambda}(x-\delta
e_j)-2\tilde{\Phi}_{\lambda}(x) \right) \longrightarrow 0~.
\end{equation}
For $x\in B$, and  $\delta>0$ such that $B\left(x,\frac{\delta}{C}\right)\subset B$, according to (\ref{Laplacien}),
$$\Delta ^{(2)}_{\delta}\tilde{\Phi}_{\lambda}(x)\le \frac{\delta^2}{2}\int _{\R^d}
(h(\xi)-1){\bf 1}_{2h(\xi)\le H}   |x-\xi|^{-d+2h(\xi)-2} ~d\xi\le 0.$$
Then (\ref{delta2phi}) implies that $\{\xi;2h(\xi)\le H\}$ has measure zero, which contradicts the assumption $H>2m$.\\
The proof of the first part of Theorem \ref{thlassh} is now complete.
\\ \\

\subsection{Fdd-lass properties for the X-ray transforms}
Finally we consider the fdd-lass property at point $y_0$ for the X-ray transform.
As previously, we restrict
the computation to the one-dimensional distribution and denote $\widetilde{{\mathcal  P}_{\alpha}X}={\mathcal  P}_{\alpha}X-\mathbf{E}\left({\mathcal  P}_{\alpha}X\right)$ the centered version of ${\mathcal  P}_{\alpha}X$. For any  $y\in <\alpha>^{\perp}$, $t\in \R$ and $H\in (0,1)$, we write
$$\mathbf{E}\exp   \left(   it\lambda^{-H}(\widetilde{{\mathcal  P}_{\alpha}X}(y_0+\lambda
  y)-\widetilde{{\mathcal  P}_{\alpha}X}(y_0)) \right)
=\exp \Phi(H,\lambda,y_0,y, t)$$
where $\Phi(H,\lambda,y_0, y, t)$ is given by
$$
\int_{\R^d\times \R^+} \left(e^{it\lambda^{-H}G_{\rho}(\lambda
y,\xi-y_0,r)}-1-it\lambda^{-H}G_{\rho}(\lambda y,\xi-y_0,r)\right) \,d\nu_h(\xi,r).$$
By the same change of variable as in the covariance part of the
proof,  $\Phi(H,\lambda,y_0, y, t)$ is equal to
\begin{eqnarray*}
&&\int_{<\alpha>^{\perp}\times \R \times (0,\lambda^{-1})} \lambda^{-1} \left(
e^{it\lambda^{1-H}G_{\rho(p+\lambda .)}(
y,\gamma,r)}-1-it\lambda^{1-H}G_{\rho(p+\lambda .)}(y,\gamma,r)\right)\\
& &~~\times ~(\lambda r)^{2h(\lambda\gamma+p\alpha+y_0)} r^{-d-1}~d\gamma dpdr~.
\end{eqnarray*}
Let us remark that
$$\lambda^{1-H}G_{\rho(p+\lambda .)} \stackrel{\lambda\rightarrow 0^+}{\longrightarrow} 0$$
and so
$$ \lambda^{-1} \left(
e^{it\lambda^{1-H}G_{\rho(p+\lambda.)}}-1-it\lambda^{1-H}G_{\rho(p+\lambda.)}
\right) \underset{\lambda\rightarrow 0^+}{\sim} -\frac{1}{2} t^2 \lambda^{-2H+1}G^2\rho^2~,$$
where $f(\lambda)\underset{\lambda\rightarrow 0^+}{\sim}g(\lambda)$ if
$\frac{f(\lambda)}{g(\lambda)} \underset{\lambda\rightarrow 0^+}{\longrightarrow}1$.
Consequently we  can argue along the  same lines as  in the covariance part  of the
proof to get, for $H=m(\alpha,y_0)+1/2$,
$$\Phi(H,\lambda,y_0,y, t)\stackrel{\lambda\rightarrow 0^+}{\longrightarrow}
-~\frac{1}{2} t^2~ \left( \int_{\R}\mathbf{1}_{h(y_0+p\alpha)=m(\alpha,y_0)}\rho(p)^2 dp\right)
\Gamma^{m(\alpha,y_0)+1/2}(y,y),$$
which concludes the proof.

\section*{Acknowlodgements}
We  would  like  to warmly  thank  Aline  Bonami  for
 her relevant
contribution to   simplifying  many computations lines, as well as for
very fruitful discussions, and Peter Scheffler for the careful reading  of our manuscript.


\end{document}